\documentclass{gtart_a}
\pdfoutput=1
\usepackage{pinlabel}

\title{Alternate Heegaard genus bounds distance}   

\author{Martin Scharlemann}
\givenname{Martin}
\surname{Scharlemann}
\address{Mathematics Department\\
        University of California\\\newline
        Santa Barbara, CA 93106\\USA}
\email{mgscharl@math.ucsb.edu}
\urladdr{}

\author{Maggy Tomova}
\givenname{Maggy}
\surname{Tomova}
\address{Mathematics Department\\
	University of Iowa\\\newline
         Iowa City, IA 52242\\USA}
\email{mtomova@math.uiowa.edu}
\urladdr{}

\volumenumber{10}
\issuenumber{}
\publicationyear{2006}
\papernumber{16}
\lognumber{0620}
\startpage{593}
\endpage{617}

\doi{}
\MR{}
\Zbl{}

\arxivreference{math.GT/0501140}
\arxivpassword{uc97r}   


\received{1 June 2005}
\revised{14 November 2005}
\published{4 May 2006}
\publishedonline{4 May 2006}
\accepted{27 March 2006}
\proposed{David Gabai}
\seconded{Joan Birman, Cameron Gordon}
\corresponding{}
\editor{}
\version{}

\arxivreference{math.GT/0501140}
\arxivpassword{uc97r}

\let\xysavmatrix\xymatrix
\def\xymatrix{\disablesubscriptcorrection\xysavmatrix}
\AtBeginDocument{\let\bar\wbar}

\makeatletter
\def\cnewtheorem#1[#2]#3{\newtheorem{#1}{#3}[section]
\expandafter\let\csname c@#1\endcsname\c@thm}

\newtheorem{thm}{Theorem}[section]  
\cnewtheorem{cor}[thm]{Corollary}
\cnewtheorem{prop}[thm]{Proposition} 
\cnewtheorem{lemma}[thm]{Lemma} 
\theoremstyle{definition}
 
\cnewtheorem{quest}[thm]{Question}
\cnewtheorem{defin}[thm]{Definition} 
\cnewtheorem{example}[thm]{Example} 
\makeatother

\newcommand{\aaa}{\mbox{$\alpha$}}
\newcommand{\map}{\mbox{$\rightarrow$}}

\newcommand{\Sss}{\mbox{$\Sigma$}}

\newcommand{\Ggg}{\mbox{$\Gamma$}}
\newcommand{\Lll}{\mbox{$\Lambda$}}
 
\newcommand{\bdd}{\mbox{$\partial$}}

\newcommand{\xX}{\mbox {\sc x}}
\newcommand{\yY}{\mbox {\sc y}} 

\begin{document}

\begin{asciiabstract} 
Suppose M is a compact orientable irreducible 3-manifold with Heegaard
splitting surfaces P and Q.  Then either Q is isotopic to a possibly
stabilized or boundary-stabilized copy of P or the distance d(P) is
less than or equal to twice genus(Q).  More generally, if P and Q are
bicompressible but weakly incompressible connected closed separating
surfaces in M then either P and Q can be well-separated or P and Q are
isotopic or d(P) is less than or equal to twice genus(Q).
\end{asciiabstract}

\begin{htmlabstract}
<p class="noindent">
Suppose M is a compact orientable irreducible 3&ndash;manifold with
Heegaard splitting surfaces P and Q.  Then either Q is isotopic to
a possibly stabilized or boundary-stabilized copy of P or the distance
d(P) &le; 2 genus(Q).</p>

<p class="noindent">
More generally, if P and Q are bicompressible but weakly
incompressible connected closed separating surfaces in M then either</p>

<ul>
<li><p class="noindent"> P and Q can be well-separated or</p>
<li><p class="noindent"> P and Q are isotopic or</p>
<li><p class="noindent"> d(P) &le; 2 genus(Q).</p>
</ul>
\end{htmlabstract}

\begin{abstract}
Suppose $M$ is a compact orientable irreducible $3$--manifold with
Heegaard splitting surfaces $P$ and $Q$.  Then either $Q$ is isotopic to
a possibly stabilized or boundary-stabilized copy of $P$ or the distance
$d(P) \leq 2\, \mathrm{genus}(Q)$.

More generally, if $P$ and $Q$ are bicompressible but weakly
incompressible connected closed separating surfaces in $M$ then either

$\bullet$\qua$P$ and $Q$ can be well-separated or\newline
$\bullet$\qua$P$ and $Q$ are isotopic or\newline
$\bullet$\qua$d(P) \leq \mathrm{genus}(Q)$.

\end{abstract}

\maketitle

\section{Introduction}

Suppose $M$ is an irreducible compact orientable $3$--manifold and $P \subset 
M$  is a closed connected separating surface
properly embedded in $M$.  $P$ is {\em bicompressible} if 
it compresses into both complementary components $A$ and $B$. $P$ is
{\em strongly compressible} if there are compressing disks for $P$ in 
$A$ and $B$ which have disjoint boundaries in $P$.  If this is not the
case, then $P$ is {\em weakly incompressible}.

Given a closed bicompressible surface $P$ in $M$, there is a natural
generalization, essentially due to Hempel \cite{He}, of the notion of
strong compressibility.  Let $\mathcal{U}, \mathcal{V}$ be the sets of
isotopy classes of essential simple closed curves in $P$ that bound
disks in $A$ and $B$ respectively.  The {\em distance} $d(P)$ is the
distance in the curve complex of $P$ between the set of vertices
corresponding to $\mathcal{U}$ and the set of vertices corresponding
to $\mathcal{V}$.  Put another way, $d(P)$ is the smallest number $n
\geq 0$ so that there is a sequence of essential simple closed curves
$\aaa_{0},\ldots, \aaa_{n}$ in $P$ with $\aaa_{0}$ bounding a disk in
$A$, $\aaa_{n}$ bounding a disk in $B$ and for each $1 \leq i \leq n$,
$\aaa_{i-1}$ and $\aaa_{i}$ can be isotoped in $P$ to be disjoint.
Note that $P$ is strongly compressible if and only if $d(P) \leq 1$.

Suppose $Q$ is another closed bicompressible surface in $M$.  Can $Q$
tell us anything about $d(P)$?  Obviously it need not: If $Q$ is
isotopic to $P$ then it can provide no new information about $d(P)$.
At the other extreme, if $Q$ can be isotoped well away from $P$ (for
example, if it is an unknotted handlebody inside a ball) then it also
can tell us nothing about $P$.  Our intention is to show that, aside
from these two possibilities,  in fact $d(P) \leq \mathrm{genus}(Q)$.  In
particular this conclusion follows when $P$ and $Q$ are non-isotopic
irreducible Heegaard splittings for $M$.  In the context of the theory
of Heegaard splittings, our central results are these:

\medskip
{\bf \fullref{cor:Heegmain}}\qua\textsl{Suppose $P$ and $Q$ are
Heegaard splitting surfaces for the compact orientable $3$--manifold
$M$.  Then either $d(P) \leq \mathrm{genus}(Q)$ or $Q$ is isotopic
to $P$ or to a stabilization or boundary-stabilization of $P$.}

\medskip
One can frame the Corollary in a different way:

{\bf Corollary}\qua {\sl If a compact orientable $3$--manifold $M$ has a 
has a genus $g$ Heegaard splitting surface $P$ with $d(P) > 2g(P)$, then
    
    \begin{itemize}
	\item $P$ is a minimal genus Heegaard splitting surface of $M$;
	\item any other Heegaard splitting surface of the same genus is isotopic to $P$; in fact,
\item any Heegaard splitting surface $Q$ with $\mathrm{genus}(Q) \leq d(P)/2$ is isotopic to a stabilization or
boundary stabilization of $P$.
	\end{itemize} }

\medskip
We would like to thank Abby Thompson for pointing out the need, when $M$ has boundary, to include the possibility that $Q$ is a boundary-stabilized copy of $P$.

\medskip

Following Kobayashi \cite{Ko}, Hempel \cite[Theorem 2.7]{He} describes
a method to construct closed $3$--manifolds that have Heegaard
splittings of any specified genus $g \geq 2$ and arbitrarily large
distance.  Then:

\medskip

{\bf Corollary}\qua {\sl For any genus $g \geq 2$ and any $k$ there is a closed
    $3$--manifold $M$ so that 
    
    \begin{itemize}
	\item $M$ has a genus $g$ Heegaard splitting $M = H_{1}
	\cup_{P} H_{2}$ and 
	\item for each $h \leq k$ any genus $h$ Heegaard splitting of
	$M$ is isotopic to $P$ or to a stabilization of $P$.
	\end{itemize} }

	\medskip
	
	Of course the second conclusion implies in particular that there are no
	Heegaard splittings of genus less than $g$, so we may as well 
	take $k \geq g$.
	
	\begin{proof} Following \cite[Theorem 2.7]{He} construct a closed
	$3$--manifold $M$ and a Heegaard splitting $M = H_{1} \cup_{P}
	H_{2}$ so that $\mathrm{genus}(P) = g$ and $d(P) > 2k$.  Suppose $Q$
	is a Heegaard splitting surface that is not isotopic to $P$ or
	to a stabilization of $P$.  Then by \fullref{cor:Heegmain}
	$2k < d(P) \leq \mathrm{genus}(Q)$ so $\mathrm{genus}(Q) > k$. 
	\end{proof}

This work was partially supported by National Science Foundation grant DMS-0405712.

\section{Preliminaries}

Throughout the paper, all $3$--manifolds will be compact and
orientable.  Many of the technical foundations for the argument here
were laid in Scharlemann \cite{Sc}, but the setting was somewhat
different.  In \cite{Sc} there was only one bicompressible surface, it
was not necessarily connected and each of its components had boundary.
In order to incorporate those results we briefly review some of the
relevant definitions and results there.

\begin{defin} \label{defin:splitting}
     Let $(Q, \bdd Q) \subset (M, \bdd M)$ be a properly embedded
    orientable surface in the irreducible $3$--manifold $M$.
    $Q$ will be called a {\em splitting surface} if no component is
    closed, no component is a disk, and $M$ is the
    union of two $3$--manifolds $X$ and $Y$ along $Q$.  
    
    We abbreviate by saying that $Q$ splits $M$ into the submanifolds
    $X$ and $Y$. \end{defin}
    
\begin{defin} \label{defin:compressiontypes}
Suppose that $(Q, \bdd Q) \subset (M, \bdd M)$ is either a splitting
surface that splits $M$ into submanifolds $X$ and $Y$ or $Q$ is a
closed connected separating surface with complementary components $X$
and $Y$.  Then $Q$ is {\em bicompressible} if both $X$ and $Y$ contain
compressing disks for $Q$ in $M$; $Q$ is {\em strongly compressible}
if there are such disks whose boundaries are disjoint in $Q$.  If $Q$
is not strongly compressible then it is {\em weakly incompressible}.
\end{defin}

Let $S_{0}, S_{1}$ be two connected compact subsurfaces in the same
component $N$ of $\bdd M$, with each component of $\bdd S_{i}, i = 0,
1$ essential in $\bdd M$ and $S_{0} \subset interior(S_{1})$.  Let
$T_{i}, i = 0, 1$ be the properly embedded surface in $M$ obtained by
pushing $S_{i}$, rel $\bdd S_{i}$, into the interior of $M$, so the
region $R$ lying between $T_{0}$ and $T_{1}$ is naturally homeomorphic
to $S_{1} \times I$. (The boundary of $R$ is the union of $T_{0},
T_{1}$ and $S_{1} - S_{0}$.) 

\begin{defin} \label{defin:tuberecess} The properly embedded surface
$T_{0} \cup T_{1} \subset M$ is called the {\em recessed collar}
determined by $S_{0} \subset S_{1}$ bounding $R$.  The properly
embedded surface obtained by tubing $T_{0}$ to $T_{1}$ along an
$I$--fiber of $S_{1} \times I \cong R$ that is incident to $T_{0}$ is
called a {\em tube-spanned recessed collar} in $M$.
\end{defin}

\begin{prop}\label{prop:essential} {\rm\cite[Proposition 2.5]{Sc}}\qua 
Suppose $M$ is an irreducible $3$--manifold, $N$ is
a compressible component of $\bdd M$ and $(Q, \bdd Q) \subset (M, \bdd
M)$ is a properly embedded essential surface containing no disk
components and with at least one essential component incident to $N$.
Let $\mathcal{V}$ be the set of essential curves in $N$ that bound
disks in $M$ and let $q$ be any component of $\bdd Q$.
 Then, either $d(\mathcal{V}, q) \leq 1-\chi(Q)$ or $q$ lies in the
    boundary of a $\bdd$--parallel annulus component of $Q$.

\end{prop}

\begin{thm}\label{thm:bicompressible} {\rm\cite[Theorem 5.4]{Sc}}\qua 
Suppose $M$ is an irreducible  $3$--manifold, $N$ is
a compressible boundary component of $M$ and $(Q, \bdd Q) \subset (M,
\bdd M)$ is a bicompressible, weakly incompressible splitting surface 
with a bicompressible component incident to $N$. 
	    
Let $\mathcal{V}$ be the set of essential curves in $N$ that bound
disks in $M$ and let $q$ be any component of $\bdd Q \cap N$.  Then either

\begin{itemize}
\item $d(q, \mathcal{V}) \leq 1 - \chi(Q)$ in the curve
	 complex on $N$ or 
\item $q$ lies in the boundary of a $\bdd$--parallel annulus component of $Q$ or 
\item one component of $Q$ is a tube-spanned recessed collar; all other
	 components of $Q$ incident to $N$ are incompressible and $\bdd$--parallel.
	 \end{itemize}   
\end{thm}

In addition, we note the following more technical lemma:

\begin{lemma} \label{lem:essentialdist}

    Suppose $S$ is a closed connected bicompressible separating
    surface in $M$, dividing $M$ into components $A$ and $B$.  Let
    $\mathcal{U}$ and $\mathcal{V}$ be the collections of simple
    closed curves in $S$ that bound disks on sides $A$ and $B$
    respectively.  Let $F \subset M$ be a properly embedded connected
    surface transverse to $S$.  Let $F_{A} = F \cap A, F_{B} = F \cap
    B$ and suppose:
   \begin{itemize}
\item All curves of $F \cap S$ are essential on $S$ and $F$.  

\item
There is at least one $a \in \bdd F_A$ such that $d(a,\mathcal{U})
\leq 1-\chi(F_A)$ and any $a \in \bdd F_A$ for which the inequality
does not hold is the boundary of an annulus component of $F_A$ which
is parallel into $S$.  

\item There is at least one $b \in \bdd F_B$
such that $d(b,\mathcal{V}) \leq 1-\chi(F_B)$ and any $b \in \bdd F_B$
for which the inequality does not hold is the boundary of an annulus
component of $F_B$ which is parallel into $S$.
   \end{itemize}

    Then $d(S) \leq 2-\chi(F)$ .

\end{lemma}
\begin{proof}
Call a component $c$ of $F \cap S$ {\em A--conforming} (resp {\em
B--conforming}) if $d(c,\mathcal{U}) \leq 1-\chi(F_A)$ (resp
$d(b,\mathcal{V}) \leq 1-\chi(F_B)$).  By hypothesis there are both
A--conforming components of $F \cap S$ and B--conforming components.  If
there is a component $c$ that is both A--conforming and B--conforming,
then $$d(S) = d(\mathcal{U}, \mathcal{V}) \leq d(c,\mathcal{U}) +
d(c,\mathcal{V}) \leq 2 - \chi(F_{A}) - \chi(F_{B}) = 2 - \chi(F)$$ as
required.  

If there is no such component, let $\gamma$ be a path in $F$ from an
A--conforming component to a B--conforming component, chosen to
intersect $S$ as little as possible.  In particular, any component of 
$F \cap S$ incident to the interior of $\gamma$ is neither A--conforming 
nor B--conforming, so each of these components of $F_{A}$ 
and $F_{B}$ is an annulus, parallel to an annulus in $S$.  It follows that the
components of $F \cap S$ at the ends of $\gamma$ are isotopic in $S$. 
Letting $c$ be a simple closed curve in that isotopy class in $S$ we
have as above $$d(S) = d(\mathcal{U}, \mathcal{V}) \leq d(c,\mathcal{U}) +
d(c,\mathcal{V}) \leq 2 - \chi(F_{A}) - \chi(F_{B}) = 2 - \chi(F)$$ as
required.  \end{proof} 

\begin{cor} \label{cor:essentialdist}

    Suppose $M$, $S$, $A$, $B$, $\mathcal{U}$ and $\mathcal{V}$ are as
    in \fullref{lem:essentialdist}.  Let $F \subset M$ be a properly
    embedded connected surface transverse to $S$ so that all curves of
    $F \cap S$ are essential on $S$ and $F$.  If $F_{A} = F \cap A$
    and $F_{B} = F \cap B$ are essential in $A$ and $B$ respectively, 
    then $d(S) \leq 2-\chi(F)$ .

\end{cor}

\begin{proof} \fullref{prop:essential} shows that $F_{A}$ and 
    $F_{B}$ satisfy respectively the second and third conditions of \fullref{lem:essentialdist}.
\end{proof}

\section{Pairs of bicompressible surfaces and their intersection
configurations}

Throughout this section, $P$ and $Q$ will be closed connected
orientable bicompressible and weakly incompressible separating
surfaces in the irreducible $\bdd$--irreducible $3$--manifold
$M$.  The complementary components of $P$ will be denoted $A$ and $B$;
those of $Q$ will be denoted $X$ and $Y$.

We introduce the following notation: Let $P^{A}$ denote the possibly
disconnected surface obtained by maximally compressing $P$ towards
$A$.  If a $2$--sphere arises and the $3$--ball it bounds in $M$ lies on
the $A$ side, delete the $2$--sphere.  That is, add the $3$--ball to the
part of $M$ that lies between $P$ and $P^{A}$.  In the end, $P^{A}$
may be empty; this is equivalent to $A$ being a handlebody.  Similarly
define $P^{B}$ by maximally compressing $P$ towards $B$, $Q^{X}$ by
maximally compressing $Q$ towards $X$ and $Q^{Y}$ by maximally
compressing $Q$ towards $Y$.

The surface $P^{A}$ separates $M$.  Let $A_{-} \subset A$ denote the
closure of the region that lies between $P$ and $P^{A}$ and similarly
define $B_{-}, X_{-}, Y_{-}$.  There is a collection of $1$--handles in
$A_{-}$, dual to the $2$--handles by which $P$ was compressed to
$P^{A}$ so that when $1$--surgery is done to $P^{A}$ along these
handles, $P$ is recovered.  Alternatively, we can think of the cores
of these $1$--handles as forming a graph $G$ in $A_{-}$, disjoint from
$P^{A}$ except precisely on its valence one vertices.  A regular
neighborhood of $G \cup P^{A}$ then has boundary isotopic to the union
of $P^{A}$ and $P$.  We will denote $G \cup P^{A}$ by $\Sss^{A}$ and
call it a {\em spine} of $A_{-}$.  There are parallel constructions of
$3$--manifolds in $M$ denoted $B_{-}, X_{-}, Y_{-}$ and their spines
$\Sss^{B}, \Sss^{X}, \Sss^{Y}$ respectively.  Notice that spines are
not well-defined as graphs, but may be altered by sliding edges along
edges or along the surface (eg $P^{A}$) to which the edge is attached.

The surfaces $P$ and $Q$ in $M$ are said to be {\em well-separated} if
they can be isotoped so that $A_{-} \cup B_{-}$ and $X_{-} \cup Y_{-}$
are disjoint.  We are interested in understanding $P \cap Q$ when
these surfaces are not well-separated, so we henceforth assume this. 
In particular, we assume that neither $P$ nor $Q$ lies in a $3$--cell, though we note below
that this possibility is ruled out anyway by the assumption that each is both
bicompressible and weakly incompressible.  

Note that if $P^{A}$ contains a sphere component, then the ball the
sphere bounds can't lie entirely in $A$ (else it would have been
removed) so in fact all of $A_{-} \cup B_{-}$ lies in the ball.  This
would mean in particular that $P$ and $Q$ are well-separated.  So
henceforth we further assume that no component of $P^{A}, P^{B},
Q^{X}$, or $Q^{Y}$ is a sphere. It follows that the only way in which
$P$ or $Q$ could be a torus is if $M$ is a Lens space (or $S^{3}$) and
$P$ and $Q$ are both genus one Heegaard surfaces for $M$.  In that
case, it is known that $P$ and $Q$ are isotopic (Bonahon--Otal
\cite{BoO}), so we further assume that neither $P$ nor $Q$ are tori.

Summarizing:

\medskip
{\sl Henceforth assume that $P$ and $Q$ are not well-separated, that
neither surface is a torus, and that no component of $P^{A}, P^{B},
Q^{X}$, or $Q^{Y}$ is a sphere.}

\medskip
Each of the surfaces $P^{A}, P^{B}, Q^{X}$, or $Q^{Y}$ is
incompressible.  For suppose, say, that $P^{A}$ or $P^{B}$ is
compressible.  Then $P^{A} \cup P^{B}$ is compressible.  The
compressing disk cannot lie in $A - A_{-}$ else $P^{A}$ could be
further compressed into $A$.  Nor can it lie in $B - B_{-}$. Hence it
lies in $A_{-} \cup B_{-}$, a manifold that is Heegaard split by $P$.
Then it is a classical result, essentially due to Haken (cf
Casson--Gordon \cite{CG}) that it is a weakly reducible splitting, ie
$P$ is strongly compressible.

Even if we allowed $P$ and $Q$ to be well-separated, this observation rules out the possibility that $P$ or $Q$ lies in a $3$--cell $C$ in $M$.
For if, say, $P \subset C$ then all compressions in the construction could be done in $C$ and, ultimately, $P^{A}, P^{B}$ would, absurdly, be incompressible surfaces in the $3$--cell $C$.

\medskip
    
Suppose that $P$ and $Q$ are in general position and there are curves
of intersection that are inessential in both surfaces.  Each such
curve $c$ bounds both a disk $D \subset P$ and a disk $E \subset Q$.
If $D$ intersects $Q$ and $E$ intersects $P$ only in curves that are
inessential in $Q$ (resp $P$) we say that $c$ is {\em removable}.  The
term reflects the fact that all such curves can be removed by
isotopies of $P$ whose support lies away from any curves of
intersection that are essential either in $P$ or in $Q$.  Indeed, if
$c$ is removable, then any component of $D \cap E$ is clearly also
removable.  An innermost one in $D$ (say) cuts off a disk $D'$ from
$D$ and a disjoint disk $E'$ from $E$; neither intersects $P$ or $Q$
in essential curves.  Since $M$ is irreducible, the union of $D'$ and
$E'$ bounds a ball in $M$, which can be used to isotope $E'$ past
$D'$, removing only removable curves.  Continue the process until all
removable curves are gone.

Having removed all removable curves of intersection, proceed to 
associate one or more of the following labels to the configuration 
given by $P$ and $Q$:

\begin{itemize}
    
    \item Label $A$ (resp $B$) if some component of $Q \cap 
    P$ is a meridian circle for $A$ (resp $B$).  
        
    \item $X$ (resp $Y$) if there is a meridian disk for $X$ (resp
    $Y$) that is disjoint from $P$.
    
    \item $x$ (resp $y$) if some spine $\Sss^{A}$ or $\Sss^{B}$ lies 
    entirely in $Y$ (resp $X)$.
    
    \end{itemize}
    
The following is immediate:

\begin{lemma} \label{lemma:nolabel} If the configuration of $P$ and
$Q$ has no label above, then $d(P) \leq 2 -
\chi(Q)$.
\end{lemma}

\begin{proof} Consider the curves $P \cap Q$ and suppose some are
essential in $P$ but inessential in $Q$.  An innermost such curve in
$Q$ will bound a disk in $A$ or $B$ (cf Robinstein--Scharlemann
\cite{RS}[Lemma 4.3]).  Since there is no label, such curves can not
exist.  In particular, any intersection curve that is inessential in
$Q$ is inessential in $P$.  Now suppose there is a curve of
intersection that is inessential in $P$.  An innermost such curve $c$
bounds a disk $D \subset P$ that lies either in $X$ or in $Y$ but,
because there is no label $X$ or $Y$, this curve must be inessential
in $Q$ as well.  Let $E$ be the disk it bounds there.  We have just
seen that all intersections of $E$ with $P$ must be inessential in
both surfaces, so $c$ is removable and would have been removed at the
onset.  We conclude that all remaining curves of intersection are
essential in both surfaces.

If $P \cap Q =\emptyset$ then $Q \subset A$ say so $B$ is entirely
contained in $X$ or in $Y$, say in $Y$.  But $B$ contains all spines
$\Sigma^B$ so there will be a label $x$ contradicting the hypothesis.
Thus $P\cap Q \neq \emptyset$.

If $Q_{A} = Q \cap A$ were compressible in $A$ the boundary of the
compression disk would also be essential in $Q$, since no curve of
$Q \cap P$ is inessential in $Q$.  In other words, there would be a
meridian disk for $X$ or $Y$ contained in $A$.  As there are no labels
$X$ or $Y$, $Q_A$ and similarly $Q_B=Q\cap B$ must be incompressible in
$A$ and $B$ respectively.  Now suppose that every component of $Q_A$
were parallel into $P$ so the region between each component of $Q_A$
and $P$ is a product.  We could then use this
product structure to isotope $Q_A$ so that it is disjoint from a spine
$\Sigma^{A}$
contradicting that the region is not labeled.  Therefore $Q_A$ and
$Q_B$ must be essential in $A$ and $B$ respectively.  The bound on the
distance then follows by \fullref{cor:essentialdist}.
\end{proof}

\begin{lemma} \label{lemma:spineinY}
If there is a spine $\Sss^{A} \subset Y$ then any meridian disk 
for $Y$ that is disjoint from $P$ must intersect $\Sss^{A}$.  
    \end{lemma}
    
\begin{proof} Since $\Sss^{A} \subset Y$, $Q$ is disjoint from the
separating surface $P^{A}$.  In particular, $Q$ either lies in $A - 
A_{-}$ or in $B \cup A_{-}$.  The former is impossible, since all 
compressions of $Q$ could then be done in $A - A_{-}$ (since $P^{A}$ 
is incompressible) implying $P$ and $Q$ are well-separated.   So $Q \subset B
\cup A_{-}$.  

Suppose $E$ is a meridian disk for $Y$ that is disjoint from $P$ and
from $\Sss^{A}$.  Use the product structure between $P$ and $\Sss^{A}$
to push all of $Q \cap A$, as well as $E$, into $B$; after that push,
maximally compress $Q \cap B$ in $B \cap Y$ to get a surface $Q'$, an
incompressible surface in $B$.  Compressing
disks for $P$ in $B$ can be taken to be disjoint from $Q'$, so $P^{B}$
can be chosen so that $Q' \cap B_{-} = \emptyset$.  That is, $B_{-}$
remains completely in $Y$ even after $Q$ is compressed sufficiently
into $Y$ to make $Q'$.  This means that $Q$ itself is disjoint from
$B_{-}$.  By the argument just given, this would imply again that $P$
and $Q$ can be well separated.
\end{proof}

Of course the symmetric statements hold if $\Sss^{A} \subset X$, 
$\Sss^{B} \subset Y$ or $\Sss^{B} \subset X$.

\begin{lemma} \label{lemma:notinproduct} Suppose there is a collar $P
\times I$ of $P$ and $Q \subset P \times I$.  Then

\begin{itemize}
    \item if $Q$ separates the ends of $P \times I$, then $P$ and $Q$ are
    isotopic and 

\item  if $Q$ fails to separate the ends of $P$ then $Q$ compresses 
to exactly one side in $P \times I$ and that side is a handlebody.

\end{itemize}
\end{lemma}
    
\begin{proof} If $Q$ is incompressible in $P \times I$ then it is
isotopic to $P$ and we have the first conclusion.
    
Suppose then that $Q$ is compressible into $Y \cap (P \times I)$. 
Maximally compress into $Y \cap (P \times I)$ and cap off any
resulting spheres with balls in $Y$.  If the result is the empty set,
then $Y$ is a handlebody completely contained in $P \times I$.  If $Q$
fails to also compress in $X \cap (P \times I)$ then we have the second
conclusion.  

If $Q$ does compress in both $Y \cap (P \times I)$ and in $X \cap (P
\times I)$ then a maximal compression in one direction or the other
won't yield a handlebody, since $P \times I$ has boundary.  So, with
no loss of generality, we are reduced to the case in which the result
of a maximal compression into $Y \cap (P \times I)$ is a possibly
disconnected but non-empty surface $Q' \subset (P \times I)$.  We 
show that this leads to a contradiction.

Let $V \subset Y$ denote the corresponding compression body lying
between $Q$ and $Q'$.  Then $V \cup X_{-} \subset M$ is a manifold
Heegaard split by $Q$ and, since $Q$ is weakly incompressible this is
a strongly irreducible splitting.  It follows from \cite{CG} that $Q'$
can't compress in this manifold; since $Q^{X}$ is incompressible in
$M$, it follows that $Q'$ doesn't compress in $V \cup X$.  Hence $Q'$
is incompressible in $P \times I$ and so consists of a number of
parallel copies of $P$.  Since $Q$ is connected and lies on one side
of $Q'$, the number of copies of $P$ in $Q'$ is one if $Q$ separates
the ends of $P \times I$ and two if it does not.

If $Q'$ were a single copy of $P$ then it would compress in $M$ on
both its sides, whereas we have established that $Q'$ is
incompressible on at least one side.  On the other hand, suppose that
$Q'$ consists of two copies of $P$ and that $Q$ lies between them. 
Then $V \cup X_{-}$ is a strongly irreducible Heegaard splitting of $P
\times I$, in which both boundary components lie in the same
compression body $V$.  These are understood \cite{STproduct}; in
particular, $X_{-}$ is the handlebody obtained by removing the
neighborhood of a single $I$--fiber of $P \times I$.  We have established
that both complementary components of $P \times I$ in $M$ lie in $Y$;
one of them is homeomorphic to $A$, the other to $B$.  In the end of
$P \times I$ (ie, one component of $Q'$) that abuts $A$, choose an
essential simple closed curve $a$ that bounds a disk $E$ in $A$, hence
in $Y - V$.  Let $c$ be an essential simple closed curve in $P$ that
intersects $a$ in at most one point.  Then the annulus $(c \times I)
\subset (P \times I)$ intersects $E$ in at most one point.  Removing a
vertical fiber from $c \times I$ turns the annulus into an essential
disk $D$ in $X$ intersecting $E$ in at most one point.  This leads to the
contradiction that $Q$ is strongly compressible.
\end{proof}

\begin{lemma} \label{lemma:notbothmixed}

  If there is a configuration labeled both $x$ and $Y$ but not $A$ or
  $B$ (or symmetrically $X$ and $y$ but not $A$ or $B$) then either
  $P$ and $Q$ are isotopic or $d(P) \leq 2 - \chi(Q)$.
     \end{lemma}

\begin{proof}  From the label $x$ we may assume, with no loss of
generality, that $\Sss^{A} \subset Y$.  From the label $Y$ we know
that $Q$ is compressible in $Y - P$.  Following \fullref{lemma:spineinY} we know that any such compression disk $E$
intersects $\Sss^{A}$.

We first argue that we may as well assume that all components of $P
\cap Q$ are essential in $P$.  For suppose not; let $c$ be the
boundary of an innermost disk $D$ in $P - Q$.  If $c$ were essential
in $Q$ then $D$ can't be in $Y$ (by \fullref{lemma:spineinY}) and so
it would have to lie in $X$.  But $D$ is disjoint from $E \subset A$,
contradicting the weak incompressibility of $Q$.  We deduce that $c$
is inessential in $Q$.  Replace the disk that it bounds in $Q$ with
$D$.  This changes the configuration, but it clearly can't remove a
label of the form $x, y, X, Y$ (though conceivably it may add such a
label as the part of $Q$ removed might intersect potential spines of
$A_{-}$ or $B_{-}$ or meridian disks for $X$ or $Y$.)  So the new
configuration satisfies the hypothesis also, and has at least one
fewer curve of intersection.  Continue until all curves of
intersection that are inessential in $P$ are removed.  Notice that the
moves on $Q$ are all isotopies, so the process does not change the
fact that $P$ and $Q$ can not be well-separated.  Suppose now that
some curve of intersection is inessential on $Q$.  By the above it
must be essential on $P$ but then an innermost such curve would give
rise to a label $A$ or $B$ contradicting the hypothesis.  Thus $Q \cap
A$ and $Q \cap B$ have no disk components and therefore have
nonpositive Euler characteristic.

By \fullref{lemma:spineinY} $Q \cap A$ does not compress into
$Y$ in $A - \Sss_{A}$ so it is either incompressible there or else
it compresses into $X$ there.  The latter would imply that $Q \cap
A$ is actually bicompressible in $A$.  We show that in either
case, for each component  $q$ of $ Q \cap P$ that is not the boundary 
of a $\bdd$--parallel annulus in $A$, $d(q, \mathcal{U}) \leq 1 - \chi(Q
\cap A)$.

{\bf Case 1}\qua $Q \cap A$ is incompressible in $A - \Sss_{A} \cong P
\times I$.

Then each component of $Q \cap A$ must be boundary-parallel into
$P$. The compression disk of $Q \cap A$ in $Y - P$ can be extended
via this parallelism to give a meridian disk for $A$ that is
disjoint from $Q \cap P$ and thus disjoint from all $q \in Q \cap
P$. Hence $d(q, \mathcal{U}) \leq 1 \leq 1 - \chi(Q \cap A)$. 

{\bf Case 2}\qua $Q \cap A$ is bicompressible in $A$.

$Q \cap A$ can't be a tube-spanned recessed collar (else it would
compress in $A - \Sss_{A}$) so, by \fullref{thm:bicompressible}, $d(q, \mathcal{U}) \leq 1 - \chi(Q \cap
A)$ for all $q \in \bdd Q$ which are not the boundary of an annulus 
parallel into $P$.

Now consider $ Q \cap B$.  It is incompressible in $B$ because a
compression into $Y$ would violate \fullref{lemma:spineinY} and a
compression into $X$ would provide a weak reduction of $Q$.  If $Q
\cap B$ is essential in $B$ then by \fullref{prop:essential}
$Q$ satisfies the hypotheses for $F$ in \fullref{cor:essentialdist} and the conclusion of that Proposition
completes the proof.  If $Q \cap B$ is not essential in $B$ then every
component of $Q \cap B$ is parallel into $P$ so $Q \subset P \times
I$.  It follows from \fullref{lemma:notinproduct} that either $Q$ is
isotopic to $P$ (and we are done) or $X$ is a handlebody contained in
$P \times I$.

Consider the latter possibility, $X$ is a handlebody contained in
a collar of $P$.  We know that some compressing disk $E$ in $Y$
lies entirely in $A$, so in particular it is disjoint from the end
of $P \times I$ that lies in $B$.  Put another way, we can push
the handlebody $X$ into $A$ along with $E$. Once both $X$ and $E$ lie 
in $A$, maximally compress $Q$ into $Y \cap A$ and call the resulting 
surface $Q'$.  $Q'$ is incompressible in $A$ (since $Q$ is weakly
incompressible), so all compressing disks for $P$ in $A$ can be made
disjoint from $Q'$ and hence from $Q$.  That is, $P^{A}$ separates $P$
from $Q$.  As usual, this would imply that $P$ and $Q$ are
well-separated, a contradiction to our underlying assumption.
\end{proof}

\begin{lemma} \label{lemma:notboth}
    
 If there is a configuration labelled both $X$ and $Y$ but not $A$ or
 $B$ then either $P$ and $Q$ are isotopic or $d(P) \leq 2 - \chi(Q)$.
    \end{lemma}
    
\begin{proof} 
Since $Q$ is weakly incompressible, any pair of meridian disks, one in
$X$ and one in $Y$, must intersect on their boundaries and so cannot
be separated by $P$.  It follows that if both labels $X$ and $Y$
appear then the only meridian disks for $X$ and $Y$ that are disjoint
from $P$ are both in $A$, say, so $Q_{A} = Q \cap A$ is bicompressible
in $A$ and $Q_{B} = Q \cap B$ is incompressible in $B$.

    Again we may as well assume that all components of $P \cap Q$ are
    essential in $P$.  For suppose not; let $c$ be the boundary of an
    innermost disk $D$ in $P - Q$.  If $c$ were essential in $Q$ then
    a disk in $B$ parallel to $D$ would provide a compression of
    $Q_{B}$ in $B$.  From this contradiction we deduce that $c$ is
    inessential in $Q$ and proceed as in the proof of \fullref{lemma:notbothmixed}. As no labels $A$ or $B$ appear, all curves 
    are also essential on $Q$ and therefor $Q \cap A$ and $Q \cap B$ have 
    nonpositive Euler characteristic.

    Following \fullref{thm:bicompressible} there are two cases to
    consider:

    {\bf Case 1}\qua $d(\mathcal{U}, q) \leq 1 - \chi(Q_{A})$ for each
    component $q$ of $P \cap Q$ that is not the boundary of an annulus
    component of $Q_{A}$ parallel into $P$.

If $Q_{B}$ is essential then this case follows much as in \fullref{lemma:notbothmixed}: by \fullref{prop:essential} $Q$
satisfies the hypotheses for $F$ in \fullref{cor:essentialdist} and the conclusion of that Proposition
completes the proof.  If $Q \cap B$ is not essential in $B$ then every
component of $Q \cap B$ is parallel into $P$ so $Q$ would be disjoint
from a spine $\Sigma^{B}$.  This would imply that $\Sigma^B$ is
contained in either $X$ or $Y$ so the region would also be labeled $x$
or $y$ and \fullref{lemma:notbothmixed} applies.

    {\bf Case 2}\qua $Q_{A}$ is the union of a tube-spanned recessed
    collar and some components parallel into $P$.

    If $Q_{A}$ is the union of a tube-spanned recessed collar and some
    some components parallel into $P$, then $\Sigma^A$ can be made disjoint
    from $Q$ as $Q_A$ is contained in a collar of $P$. Thus \fullref{lemma:notbothmixed} applies again.
\end{proof}

\begin{lemma} \label{lemma:notbothsmall}
  If there is a configuration labeled both $x$ and $y$ but not $A$ 
  or $B$, then either $P$
  and $Q$ are isotopic or $d(P) \leq 2 - \chi(Q)$.
     \end{lemma}

\begin{proof}
    As usual, we can assume that all curves in $P \cap Q$ are
    essential in both surfaces.  Indeed, if there is a curve of
    intersection that is innessential in $P$ then an innermost one
    either is inessential also in $Q$, and can be removed as described
    above, or is essential in $Q$ and so would give a rise to a label
    $X$ or $Y$.  In the latter case, the result would follow from
    \fullref{lemma:notbothmixed}.  As no labels $A$ or $B$ appear,
    we can again assume that all curves are also essential on $Q$ and
    thus $Q_A$ and $Q_B$ have no disk components.
    
    Both $X$ and $Y$ contain entire spines of $A_{-}$ or
$B_{-}$, though since we are not dealing with fixed spines the labels
could arise if there are two distinct spines of $A_{-}$, say, one in
$X$ and one in $Y$.  Indeed that is the case to focus on, since if
spines $\Sss_{A} \subset X$ and $\Sss_{B} \subset Y$ then $Q$ lies
entirely in $P \times I$ and separates its ends.  By \fullref{lemma:notinproduct} this implies that $Q$ is isotopic to $P$.

So suppose that $\Sss_{A} \subset Y$ and another $\Sss'_{A} \subset
X$.  Suppose $P^{A} \neq \emptyset$ (ie, $A$ is not a handlebody).
Then the parallel incompressible surfaces ${P'}^{A}$ and $P^{A}$ lie
respectively in $X$ and $Y$.  That is, $Q$ lies in a collar $P^{A}
\times I$.  Since $P^{A}$ is incompressible in $M$, $Q$ can be fully
compressed into both $X$ and $Y$ entirely within $P^{A} \times I$.
Hence $Q^{X}$ is parallel to ${P'}^{A}$, $Q^{Y}$ is parallel to
$P^{A}$ and $Q$ is a Heegaard splitting of the product $P^{A} \times
I$ and separates its ends.  But by Scharlemann--Thompson
\cite{STproduct} this would imply that $Q$ is strongly compressible, a
contradiction.

So we conclude that $A$ is a
handlebody.  If $Q \cap A$ compresses in $X \cap A$ (resp $Y \cap A$)
then the configuration has both labels $X$ and $y$ (resp $Y$ and $x$),
and we are done via \fullref{lemma:notbothmixed}.  So we henceforth
assume that $Q \cap A$ is incompressible in $A$.  Then it is certainly
incompressible in the product $A - \Sss_{A}$ and so every component of
$Q \cap A$ is parallel in $A - \Sss_{A}$ to a subsurface of $P$.
Similarly every component of $Q \cap A$ is parallel in $A - \Sss'_{A}$
to a subsurface of $P$.

Any arc in $A$ between $\Sss_{A} \subset Y$ and $\Sss'_{A} \subset X$
must intersect $Q \cap A$ an odd number of times, so there is a
component $Q_{0}$ of $Q \cap A$ that such an arc intersects an odd
number of times.  Put another way: $Q_{0}$ is a component of $Q \cap
A$ that lies between $\Sss_{A}$ and $\Sss'_{A}$.  This implies that
$Q_{0}$ is parallel into $P$ on both its sides, ie that $A
\cong Q_{0} \times I$.  Since $P$ is not a torus, $Q_{0}$ is not an
annulus.

Let $\aaa$ be an essential arc in $Q_{0}$.  Then $\aaa \times I
\subset Q_{0} \times I \cong A$ is a meridian disk $D$ for $A$ that
intersects $Q_{0}$ precisely in $\aaa$.  $\bdd$--compressing $Q_{0}$
along one of the two disk components of $D - \aaa$ gives a curve or
curves in $Q$ that are essential in $Q$ (since $Q_{0}$ is not an
annulus) and are disjoint both from $\bdd Q_{0}$ and the meridian $D
\cup_{\aaa} D'$.  In particular, for $q_{0}$ any $\bdd$--component of
$Q_{0}$, $d(q_0, \mathcal{U}) \leq 2 \leq 1 - \chi(Q_{0})$.  If
$Q_A-Q_0$ has a non-annulus component, then $1-\chi(Q_0) <
1-\chi(Q_A)$ and thus for any component $q$ in $P \cap Q$,
$d(\mathcal{U}, q) \leq d(\mathcal{U}, q_0)+d(q,q_0)\leq 1-\chi(Q_0)+1
\leq1-\chi(Q_A)$.  Thus $Q_A$ always satisifies the hypothesis for
$F_{A}$ in \fullref{lem:essentialdist}.

Now consider $Q_{B} = Q \cap B$.  If $Q_{B}$ is essential, then by
\fullref{prop:essential} $Q_{B}$ satisfies the hypothesis for
$F_{B}$ in \fullref{lem:essentialdist} and we are done by that
lemma. If $Q_{B}$ is compressible in $B$,
we have labels $X$ and $y$ (or $x$ and $Y$) and we are done via \fullref{lemma:notbothmixed}.  Finally, if each component of $Q_{B}$ is
parallel to a subsurface of $P$, then $Q$ is disjoint from a spine
$\Sigma^{B}$ as well, a case we have already considered.
\end{proof}

\section{How labels change under isotopy}

Consider how configurations and their labels change as $P$ is isotoped
in $Q$.  Clearly if there are no tangencies of $P$ with $Q$ during the
isotopy then the curves $P \cap Q$ change only by isotopies in $P$ and
$Q$, so there is no change in labels.  Similarly, if there is an index
$0$ tangency, $P \cap Q$ changes only by the addition or deletion of a
removable curve.  Since all such curves are removed before labels are 
defined, again there is no affect on the labelling.  So it suffices to
consider only what can happen to the labelling when passing through a 
saddle tangency of $P$ with $Q$.

\begin{figure}[ht!]\small
\centering
\labellist
\pinlabel $A$ [br] <0pt, -2pt> at 343 41
\pinlabel $B$ [tr] <0pt, 2pt> at 256 42
\pinlabel $C$ [b] <3pt,0pt> at 207 107
\pinlabel $P$ at 67 82
\pinlabel $Q$ at 151 15
\pinlabel $X$ at 395 95
\pinlabel $Y$ at 395 58
\endlabellist
\includegraphics[width=0.7\textwidth]{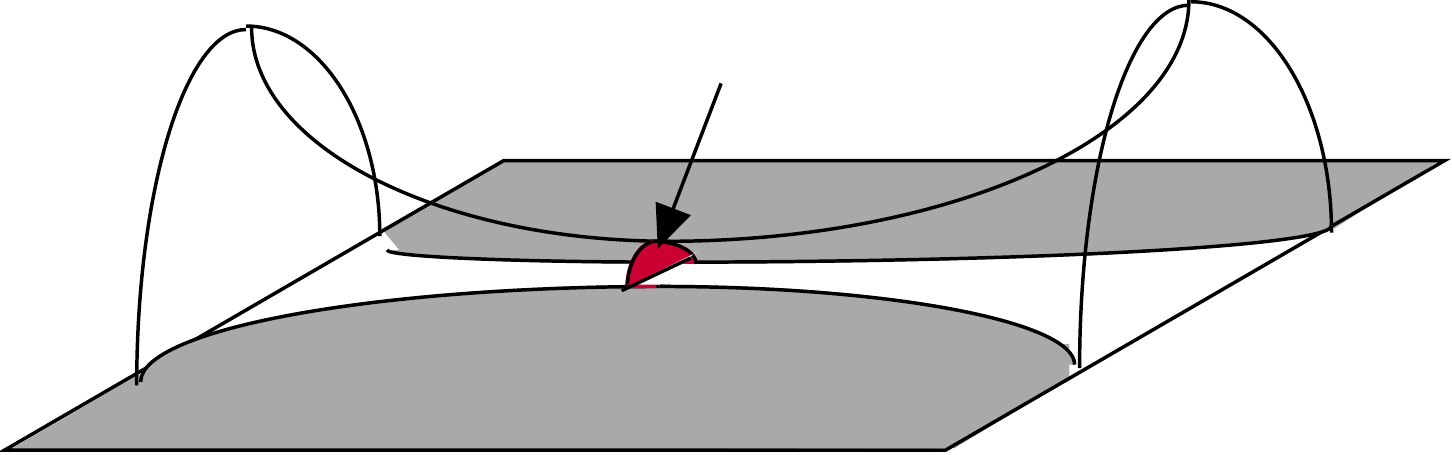}
\caption{} \label{fig:saddle}
\end{figure}

\begin{lemma} \label{lemma:saddle}
 Suppose a configuration is changed by passing through a saddle point,
 and the bigon $C$ defining the saddle tangency (cf \fullref{fig:saddle})  lies
 in $X \cap A$.  Then
 
 \begin{itemize}
 \item  No label $x$ or $X$ is removed.
 
 \item  No label $y$ or $Y$ is created.
 
 \item Suppose that, among the labels both before and after the move,
 neither $A$ nor $B$ occur.  If there is no label $x$ or $X$ before
 the move, but one is created after and if there is a label $y$ or $Y$
 before the move and none after, then either $P$ and $Q$ are isotopic
 or $d(P) \leq 2 - \chi(Q)$.
 
 \end{itemize}
 
     \end{lemma}
    
\begin{proof} We first show that no label $x$ or $X$ is removed.  If 
there is a meridian disk for $X$ that lies in $A$, a standard 
innermost disk, outermost arc argument on its intersection with $C$ 
shows that there is a meridian disk for $X$ in $A$ that is disjoint 
from $C$.  The saddle move has no effect on such a disk (nor, clearly, 
on a meridian disk for $X$ that lies in $B$).  If there is a spine of 
$A_{-}$ or $B_{-}$ lying entirely in $Y$ then that spine, too, is
unaffected by the saddle move.
    
Dually, no label $y$ or $Y$ is created: the inverse saddle move, 
restoring the original configuration, is via a bigon that lies in $B 
\cap Y$.
    
The third item is the most difficult to prove.  A useful approach is
to imagine positioning $Q$ so that it is exactly tangent to $P$ at the
saddle.  A bicollar of $Q$ then has ends which correspond to the
position of $Q$ just before the move and just after.  Let $Q_{a}$
denote $Q \cap A$ after the move and $Q_{b}$ denote $Q \cap B$ before
the move.  The bicollar description shows that $Q_{a}$ and $Q_{b}$
have disjoint boundaries in $P$.  Moreover the complement of $Q_{a}
\cup Q_{b}$ in $Q$ is a regular neighborhood of the singular component
of $P \cap Q$, with Euler characteristic $-1$.  It follows that
$\chi(Q_{a}) + \chi(Q_{b}) = \chi(Q) + 1.$

With $Q$ positioned as described, tangent to $P$ at the saddle point 
but otherwise in general position, consider the closed (non-singular) 
curves of intersection.  

{\bf Claim 1}\qua It suffices to consider the case in which all non-singular 
curves of intersection are essential in $P$.  

To prove the claim, suppose a non-singular curve is inessential and
consider an innermost one.  Assume first that the disk $D$ that it
bounds in $P$ does not contain the singular curve $s$ (ie the
component of $P \cap Q$, homeomorphic to a figure 8, that contains the
saddle point).  If $\bdd D$ is essential in $Q$, then it would give
rise to a label $X$ or a label $Y$ that persists from before the move
until after the move, contradicting the hypothesis.  Suppose on the
other hand that $\bdd D$ is inessential in $Q$ and so bounds a disk $E
\subset Q$.  All curves of intersection of $E$ with $P$ must be
inessential in $P$, since there is no label $A$ or $B$.  It follows
that $\bdd D = \bdd E$ is a removable component of intersection so the
disk swap that replaces $E$ with a copy of $D$, removing the curve of
intersection (and perhaps more such curves) has no effect on the
labelling of the configuration before or after the isotopy. So the
original hypotheses are still satisfied for this new configuration of 
$P$ and $Q$.  

Suppose, on the other hand, that an innermost non-singular innesential
curve in $P$ bounds a disk $D$ containing the singular component $s$.
When the saddle is pushed through, the number of components in $s$
switches from one $s_{0}$ to two $s_{\pm}$ or vice versa.  All three
curves are innesential in $P$ since they lie in $D$; $s_{0}$ and at
least one of $s_{\pm}$ bounds a disk in $P$ that is disjoint from $Q$.
If either is essential in $Q$ then they determine a label $X$ or $Y$
that exists both before and after the move across the saddle (since in
$Q$ the curves can be isotoped away from $s$).  This would violate the
hypotheses of the third case above, so we can assume both curves,
hence all three curves, are inessential in $Q$.  As above, this
implies that they are all removable, so passing through the
singularity can have no effect at all on the labelling.  This proves
the claim.

{\bf Claim 2}\qua It suffices to consider the case in which also all 
three curves $s_{0}, s_{\pm}$ are essential in $P$.
    
The case in which all three curves are inessential in $P$ is covered
in the proof of Claim 1.  If two are inessential in $P$, so is the
third.  We're left with the case in which exactly one is inessential
in $P$ and, following Claim 1, the disk it bounds in $P$ is disjoint
from $Q$.  If the curve were essential in $Q$ then there would have to
be a label $X$ or $Y$ that occurs both before and after the saddle
move, a contradiction.  If the curve is inessential in $Q$ then it is
removable.  If this removable curve is $s_{\pm}$ then passing through
the saddle can have no effect on the labelling; if this removable
curve is $s_{0}$ then the curves $s_{\pm}$ are parallel in both $P$
and $Q$.  In the latter case, passing through the saddle has the same
effect on labelling as passing an annulus component of $P_{X}$ across
a parallel annulus component $Q_{0}$ of $Q_{A}$.  This move can have
no effect on labels $x$ or $y$.  A meridian disk $E$ for $Y$ that's
disjoint from $P$ would persist after this move, unless $\bdd E$ is in
fact the core curve of the annulus $Q_{0}$.  But then the union of $E$
and half of $Q_{0}$ would be a meridian disk of $A$ bounded by a
component of $\bdd Q_{0} \subset P$.  In other words, there would have
to have been a label $A$ before the move, a final contradiction
establishing Claim 2.

\medskip

Claims 1 and 2, together with the fact that neither labels $A$ nor $B$
appear, reduce us to the case in which all curves of intersection are
essential in both surfaces both before and after the saddle move.
Note then that $Q_{a}$ and $Q_{b}$ are incompressible in $A$ and $B$
respectively.  For example, if the latter compresses in $B$ the
compression must be into $Y$ (since no label $X$ exists before the
move) and such a compression would persist (and so then would the
label $Y$) after the move.  Then also neither $Q_{a}$ nor $Q_{b}$ can
be inessential in $A$ or $B$ respectively.  For example, if
all components of $Q_{b}$ are parallel into $P$ then $Q_{b}$ is
disjoint from some spine of $B$ and such a spine will be unaffected by
the move, resulting on the same label ($x$ or $y$) arising before and
after the move.  We deduce that $Q_{a}$ and $Q_{b}$ are essential
surfaces in $A$ and $B$ respectively.  

Now apply \fullref{prop:essential} to both sides: Let $q_{a}$
(resp $q_{b}$) be a boundary component of an essential component of
$Q_{a}$ (resp $Q_{b}$).  Then $$d(P) = d(\mathcal{U}, \mathcal{V}) \leq
d(q_{a}, \mathcal{U}) + d(q_{a}, q_{b}) + d(q_{b}, \mathcal{V}) $$ 
$$\leq 3 - \chi(Q_{a}) - \chi(Q_{b}) = 2 - \chi(Q)$$ as required.
    \end{proof}
    
\begin{cor} \label{corollary:notboth} If two configurations are
related by a single saddle move and the union of all labels for both
configurations contains both $x$ or $X$ and $y$ or $Y$, but neither
$A$ nor $B$ then either $P$ and $Q$ are isotopic or $d(P) \leq 2 -
\chi(Q)$
    \end{cor}
    
\begin{proof} With no loss, the saddle move is as described in \fullref{lemma:saddle}.  That lemma shows that either we are done, or
there is a single configuration for which both $x$ or $X$ and $y$ or
$Y$ appear.  The result then follows from one of Lemmas
\ref{lemma:notboth}, \ref{lemma:notbothsmall} or
\ref{lemma:notbothmixed}
    \end{proof}
    
Here is our main theorem:

\begin{thm} \label{theorem:distance} Suppose $P$ and $Q$ are
bicompressible but weakly incompressible closed connected surfaces in
$M$.  Then either

\begin{itemize}
	
	\item $P$ and $Q$ can be well-separated,
	
	\item $P$ and $Q$ are isotopic or
	
	\item $d(P) \leq 2 - \chi(Q)$.
	
	\end{itemize}
	
\end{thm}

\begin{proof} We use the notation above, eg, $P$ splits $M$ into 
$A$ and $B$ containing spines $\Sss^{A}$ and $\Sss^{B}$ 
respectively.  We will assume that $P$ and $Q$ are not well-separated.

Consider a square $I \times I$ that describes generic sweep-outs of $P$ and $Q$
from $\Sss^{A}$ to $\Sss^{B}$ and from $\Sss^{X}$ to $\Sss^{Y}$
respectively.  We picture the limit of the sweep-out as 
$P$ gets near $\Sss^{A}$ as represented by the left side of the 
square.  Similarly $\Sss^{B}$ is represented by the right side of the 
square, $\Sss^{X}$ as represented by the bottom side of $I \times I$ 
and $\Sss^{Y}$ as represented by the top.  See \fullref{fig:graphic}; what's in the interior of the square in that figure
will now be explained.

Each point in the square represents a positioning of $P$ and $Q$.
Inside the square is a graph $\Ggg$, called the {\em graphic} that
represents points at which the intersection is not generic: There is a
single point of tangency between $P$ and $Q$ at each point in an edge
in the graphic; at each (valence four) vertex of $\Ggg$ there are
two points of tangency.  There are also valence two vertices
corresponding to ``birth-death'' points, but they play no role in the
discussion and we ignore them.  By general position of, say, the spine
$\Sigma^{A}$ with the surface $Q$ the graphic $\Ggg$ is incident to
$\bdd I \times I$ in only a finite number of points (corresponding to
tangencies between eg $\Sigma^{A}$ and $Q$).  Each such point in $\bdd
I \times I$ is incident to a single edge of $\Ggg$.

Any point in the complement of $\Ggg$ represents a generic 
intersection of $P$ and $Q$.  Each component of the graphic 
complement will be called a {\em region}; any two points in the same 
region represent isotopic configurations, so we can label each region 
with labels as described above. See \fullref{fig:graphic}.
If any region is unlabelled we are 
done, by \fullref{lemma:nolabel}.  

\begin{figure}[ht!]\small
\centering
\labellist
\pinlabel $\Sigma_X$ [l] <2pt, -1.5pt> at 235 66
\pinlabel $\Sigma_Y$ [l] at 235 238
\pinlabel $\Sigma_A$ [t] <3pt,0pt> at 17 26
\pinlabel $\Sigma_B$ [t] <3pt,0pt> at 199 26
\pinlabel $y$ at 175 226
\pinlabel $y$ at 28 228
\pinlabel $x$ at 32 70
\pinlabel $x$ at 185 75
\tiny
\pinlabel $A{,}Y$ at 57 231
\pinlabel $A{,}Y$ at 55 190
\pinlabel $B{,}Y$ at 143 233
\pinlabel $B{,}X$ at 130 123
\pinlabel $B{,}X$ at 171 150
\pinlabel $A{,}X$ at 87 159
\pinlabel $A{,}X$ at 73 133
\pinlabel $A{,}X$ at 44 108
\pinlabel $A$ at 43 203
\pinlabel $A$ at 32 184
\pinlabel $A$ at 50 158
\pinlabel $B$ at 100 200
\pinlabel $B$ at 180 187
\pinlabel $B$ at 193 167
\pinlabel $B$ at 164 126
\pinlabel $B$ at 180 107
\pinlabel $B$ at 131 84
\pinlabel $Y$ at 104 232
\pinlabel $Y$ at 147 196
\pinlabel $Y$ at 139 189
\pinlabel $Y$ at 161 180
\pinlabel $X$ at 62 70
\pinlabel $X$ at 92 73
\pinlabel $X$ at 94 101
\pinlabel $X$ at 129 68
\pinlabel $X$ at 148 76
\endlabellist
\includegraphics[width=0.6\textwidth]{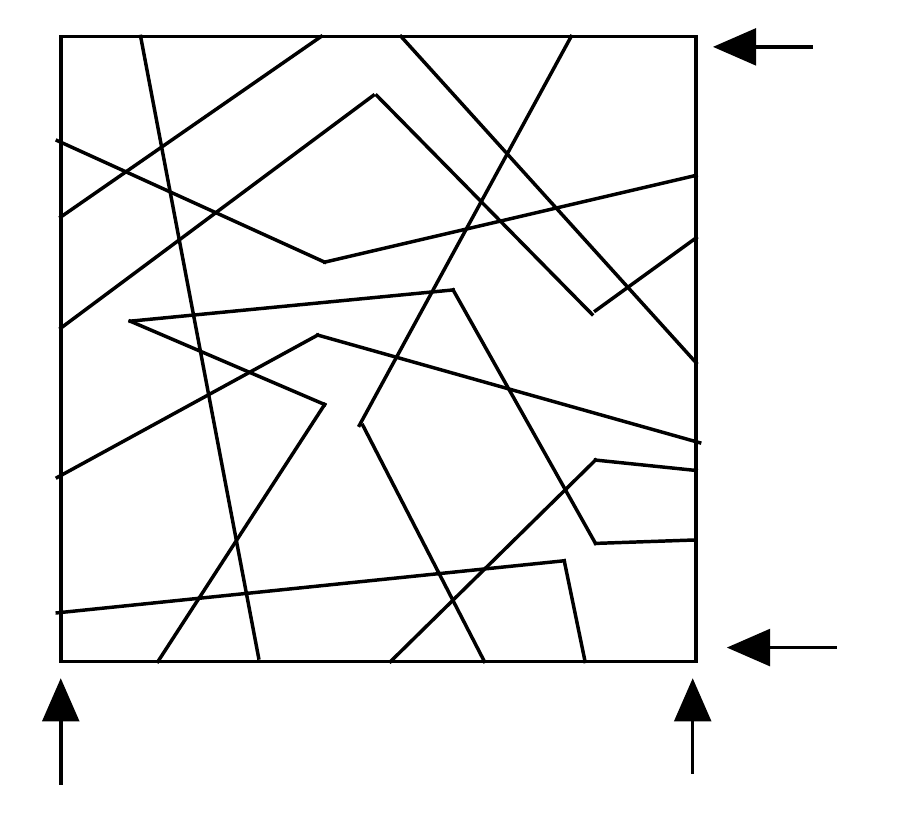}
\caption{} \label{fig:graphic}
\end{figure}  

Consider what we know about the labelling near the edges of the
square, that is when $P$ or $Q$ are near one of the spines.  For
any generic $Q$, $\Sss^{A}$ is transverse to $Q$, so $Q$ is either
disjoint from the edges in $\Sss^{A}$ or intersects them
transversally.  On the other hand, any curve in $P$ that bounds a
meridian of $B$ will intersect every meridian of $A$ (since $P$ is
weakly incompressible) so as $P$ limits on $\Sss_{A}$, any meridian
curve of $B$ will limit on a curve that runs multiple times along each
edge in $\Sss^{A}$.  It follows that when $P$ is near $\Sss^{A}$, and
for a generic $Q$ there may or may not be a label $A$, depending on
how $Q$ intersects $\Sss^{A}$.  But note that by general position all
but a finite number of meridians of edges of $\Sss^{A}$ are disjoint
from $Q$.  Hence in any case there is no label $B$, for such a label
requires a curve of $P \cap Q$ to cross all meridians of edges of
$\Sss^{A}$.

Dually, when $P$ is near $\Sigma^{B}$ there is no label $A$.  

Suppose next that $Q$ is near $\Sss^{X}$.  All but a finite number of
points on the edges of $\Sss^{X}$ will be disjoint from $P$; it
follows that for $Q$ near $\Sss^{X}$, there are meridians of $X$ that
are disjoint from $P$.  Hence all regions on the side of the square
corresponding to $Q$ near $\Sss^{X}$ (the bottom in the Figure) have a
label $X$, as well as possibly other labels.  Similarly all regions on
the opposite side of the square (the top) corresponding to $Q$ near
$\Sss^{Y}$ have at least label $Y$.

Now we proceed in almost the same way as in \cite[Section 5]{RS}.  Let
$\Lll$ be the dual complex of $\Ggg$ in $I \times I$.  Specifically,
$\Lll$ is constructed as follows: place one vertex in the interior of
each face and in the interior of each edge, including edges along
$\bdd I \times I$ (ie components of $\bdd I \times I - \Ggg$).  Then
in each face connect by an edge the vertex in the interior of the 
face with the vertex in each edge of the boundary of the face.  Each 
vertex of $\Lll$ lying on an interior edge $e$ of $\Ggg$ will then be 
incident to exactly two edges, corresponding to the faces of $\Ggg$ to 
which $e$ is incident.  Amalgamate these two edges, removing the 
vertex. Then in the end, $\Lll$ has one vertex in each face of 
$\Ggg$ and one vertex in each component of $\bdd I \times I - \Ggg$. 
Each edge of $\Ggg$ not incident to $\bdd I \times I$ crosses exactly 
one interior edge of $\Ggg$.  Each component of $I \times I - \Lll$ is 
called a face of $\Lll$.  Note that because each internal vertex 
of $\Ggg$ has valence four, each face of $\Lll$ is a quadrilateral.  
That is, $\Lll$ can be thought of as a tiling of $I \times I$ by 
quadrilaterals, where any two quadrilaterals are either disjoint, 
meet in exactly a face of each, or in a corner of each.  Arbitrarily 
many tiles may meet at a single corner.  See \fullref{fig:dual}.

\begin{figure}[ht!]
\centering
\includegraphics[width=0.4\textwidth]{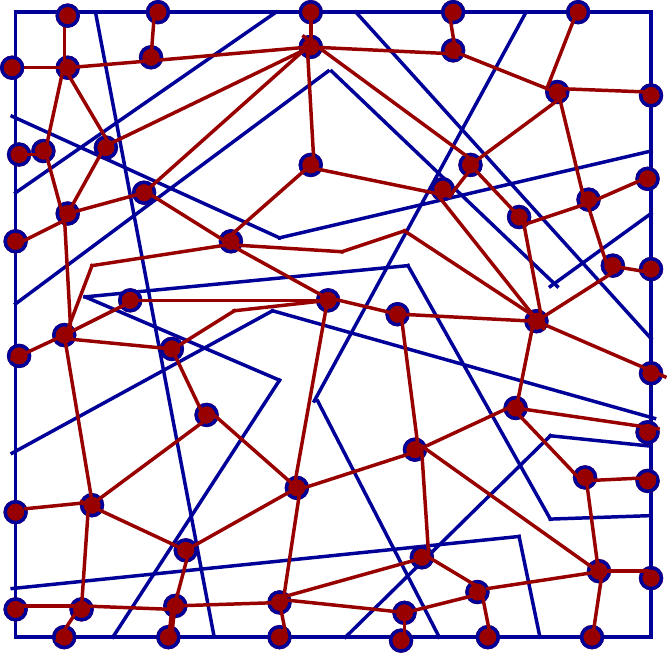}
\caption{} \label{fig:dual}
\end{figure}

Now choose the following labelling scheme for each vertex in $\Lll$:
The vertex lies in a region $R$ of $\Ggg$ which has a collection of
labels.  If a label $A$ or $B$ appears among the labels of $R$,
associate only those labels to the vertex.  If no label $A$ or $B$
appears among the labels of $R$, then associate the label $\xX$ to the
vertex if $x$ or $X$ appears among the labels of $R$ and the label
$\yY$ to the vertex if $y$ or $Y$ appears among the labels of $R$.

The first observation is that unless $P$ is isotopic to $Q$ or $d(P)
\leq 2 - \chi(Q)$ (and we are done), each vertex in $\Lll$ has exactly
one label $A, B, \xX$ or $\yY$.  Indeed, \fullref{lemma:nolabel}
says that every vertex has some label, and Lemmas
\ref{lemma:notbothmixed}, \ref{lemma:notboth} and
\ref{lemma:notbothsmall} say that no vertex of $\Lll$ has both a label
$\xX$ and $\yY$. No region can have both labels $A$ and $B$ else $P$ 
would be strongly compressible, and no region can have both a label 
$A$ or $B$ and a label $\xX$ or $\yY$ by our choice of labelling scheme.

Moreover, \fullref{corollary:notboth} says that vertices
labelled $\xX$ and $\yY$ can't even be connected by an edge in $\Lll$.
Vertices labelled $A$ and $B$ also cannot be connected by an edge in 
$\Lll$ else $P$ would be strongly compressible.  Finally, the 
observations above about labelling along the edges of $I \times I$ 
shows that no label $B$ appears along the $\Sigma^{A}$ side of $I \times 
I$ (the left side in the figure), no label $A$ appears along the 
$\Sigma^{B}$ side (the right side), no label $\yY$ appears along 
the $\Sigma^{X}$ side (the bottom) and no label $\xX$ appears along 
the $\Sigma^{Y}$ side (the top).  

We now appeal to the following quadrilateral variant of
Sperner's Lemma:

\begin{lemma} \label{lemma:sperner}  Suppose a square $I \times I$ is tiled by quadrilaterals
so that any two that are incident meet either in a corner of each or
in an entire side of each.  Let $\Lll$ denote the graph in $I \times
I$ that is the union of all edges of the quadrilaterals.  Suppose each
vertex of $\Lll$ is labelled $N, E, S,$ or $W$ in such a way that

\begin{itemize}  
    
    \item  no vertex on the East side of $I \times I$ is labelled $W$, no vertex
on the West side is labelled $E$, no vertex on the South side is
labelled $N$ and no vertex on the North side is labelled $S$.
    \item no edge in $\Lll$ has ends labelled $E$ and $W$ nor 
    ends labelled $N$ and $S$.
\end{itemize}

Then some quadrilateral contains all four labels

\end{lemma}

Two different proofs are given in a brief appendix below.

\begin{figure}[ht!]\small
\centering
\labellist
\pinlabel $X$ at 71 39
\pinlabel $Y$ at 71 109
\pinlabel $A$ at 39 74
\pinlabel $B$ at 108 74
\endlabellist
\includegraphics[width=0.3\textwidth]{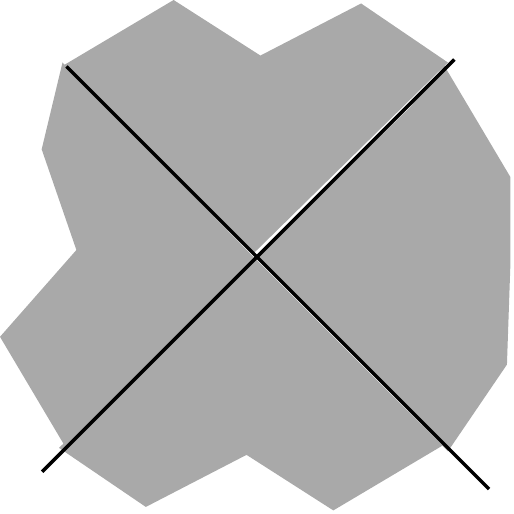}
\caption{} \label{fig:vertex}
\end{figure}

In our context the lemma says that there are four regions in 
the graphic incident to the same vertex of $\Ggg$ labelled $A, B, 
\xX$ and $\yY$.  (See \fullref{fig:vertex}).  Note then that only two saddle moves are needed to
move from a configuration labelled $A$ to one labelled $B$.  The
former configuration includes a meridian circle of $A$ and the latter
a meridian of $B$.  If, in either the region labeled $\xX$ or the region labeled $\yY$, there is a curve of intersection that is essential in $P$, then the two saddle moves show that $d(P) \leq 2 \leq 2 - \chi(Q)$, as required. 

Suppose then that in both regions labeled $\xX$ and $\yY$ all curves
of intersection are inessential in $P$.  It is easy to see (cf the
discussion before \cite{RS}[Lemma 5.6]) that at the vertex of $\Ggg$
the two saddle tangencies lie on the same singular component of $P
\cap Q$.  In particular, corresponding to the region labeled $\xX$
there is a curve in $P \cap Q$, inessential in $P$, to which separate
bands are attached, one giving rise to the label $A$ and one giving
rise to the label $B$.  (Note that if one of these bands connected
{\em different} inessential curves, the resulting curve would still be
inessential in $P$ and so could not give rise to either label $A$ or
$B$.)  Our assumption is that attaching both labels simultaneously,
corresponding to the region labeled $\yY$, gives rise to only
inessential curves in $P$.  The only way to do two simultaneous
essential band moves to an inessential curve and recover only
inessential curves is if the union of the two bands is the spine of a
torus, ie $P$ must be a torus, contradicting our initial assumption.
\end{proof}

\begin{cor} \label{cor:Heegmain} Suppose $P$ and $Q$ are both Heegaard splitting
   surfaces for 
the closed orientable $3$--manifold $M$.  Then either $d(P) \leq 2 
\mathrm{genus}(Q)$ or $Q$ is isotopic to $P$ or to a stabilization of $P$.
\end{cor}

\begin{proof} If $P$ is weakly reducible then $d(P) \leq 1$ and we are
done.  Since any Heegaard splitting of a reducible $3$--manifold is
reducible, we may as well assume that $M$ is irreducible and $P$ is
strongly irreducible.  We may as well assume $Q$ is not stabilized, so
if $Q$ is weakly reducible then there is an incompressible surface $Q'$ in
$M$ with smaller genus \cite{CG}.  In this case, the result follows from
Hartshorn's theorem \cite{Ha} (see also \cite{Sc}).  So we may as
well assume that  also $Q$ is strongly irreducible.

$P$ and $Q$ are not well-separated, for if they were there would be an
incompressible closed surface lying entirely
on one side of $Q$, ie inside a handlebody.  Since $P$ and $Q$ are strongly irreducible Heegaard surfaces (thus weakly
incompressible surfaces) and are not well-separated, the result follows
immediately from \fullref{theorem:distance}.
\end{proof}

Much the same is true when $M$ has boundary, but high distance in $P$
is also possible if $Q$ is a boundary-stabilization of $P$ (cf Moriah
\cite{Mo}).  We briefly explain what this means: Suppose $M = A \cup_P
B$ is a Heegaard splitting for $M$ and $N$ is a boundary component of
$M$ lying in $A$.  Then, since $A$ is a compression body, there is a
collar $N \times I$ embedded in $A$ so that the collar is one of the
components that remain when $P$ is maximally compressed into $A$.  The
end $N \times \{ 0 \}$ is identified with $N \subset \bdd M$ and the
complement $A - (N \times I)$ is incident to the collar via
$1$--handles attached at $N \times \{ 1 \}$.  Consider a spanning arc
$\{ * \} \times I \subset N \times I$ that has one end on $P$ and the
other end on $N \subset \bdd M$.  Add a neighborhood of $* \times
[1/2, 1]$ to the compression body $B$ and also add the subcollar $N
\times [0, 1/2]$.  Symmetrically, delete the same from $A$.  Call the
resulting $3$--manifolds $A'$ and $B'$ respectively.

The resulting decomposition $M = A' \cup_{P'} B'$ is also a Heegaard
splitting of $M$, called the {\em boundary stabilization} of $A \cup_P
B$ along $N$.  To see this, note that $A'$ is still a compression body
because the complement $N - \eta(*)$ of a disk in a closed surface is
a $2$--dimensional handlebody, so $(N - \eta(*)) \times [1/2, 1]$ is a
$3$--dimensional handlebody.

\medskip

The following Lemma, a version of the central theorem of \cite{CG} for the case of $3$--manifolds with boundary, was pointed out to us by Abby Thompson:

\begin{lemma} \label{lemma:boundstab}  Suppose $M = X \cup_Q Y$ is a genus $g$ Heegaard splitting of a $3$--manifold $M$  Then either

\begin{itemize}
 \item $Q$ is strongly irreducible
 \item $Q$ is stabilized
 \item $Q$ is boundary stabilized
 \item $M$ contains an essential closed surface of genus $< g$.
\end{itemize}
\end{lemma}

\begin{proof}  The argument is essentially in \cite{Mo}: 
Unless $Q$ is strongly irreducible, maximally weakly reduce the surface,
producing the surface $F$.  It is shown in \cite{CG} that, unless $Q$ is stabilized, $F$ is
incompressible.   The result of the disjoint
compressions of $Q$ into $X$ and $Y$ are disks in $F$ which we imagine colored
respectively red and blue.  Since tubes compressed in $Q$ only connect disks
of the same color, and $Q$ is connected, some component $F_0$ of $F$ has both
red and blue disks. If $F_0$ is not boundary parallel then it is an
essential closed surface of genus $< g$.  So suppose $F_0$ is parallel to a
boundary component via a collar $F_0 \times I$.  The original splitting is
obtained by amalgamating along $F_0$ a splitting $Q'$ for $M$ and a splitting on
$F_0  \times I$.  Since both colors of disks appear on $F_0$ the splitting on $F_0 \times
I$ is non-trivial.  According to the classification of splittings of $F_0 \times
I$ (cf \cite{STproduct}) the only non-stabilized non-trivial splitting on $F_0  \times I$ is obtained by
a single vertical tubing of collars of both its ends. Amalgamating this
splitting corresponds to boundary stabilization of $Q'$.  In other words, $Q$ is the
boundary stabilization of $Q'$ along a boundary component.\end{proof}

\begin{cor} \label{cor:Heegmain2} Suppose $P$ and $Q$ are both Heegaard splitting
   surfaces for 
the compact orientable $3$--manifold $M$.  Then either $d(P) \leq 2 
\mathrm{genus}(Q)$ or $Q$ is isotopic to $P$ or to a stabilization or boundary-stabilization of $P$.
\end{cor}

\begin{proof} Just as in the proof of \fullref{cor:Heegmain} we may as well assume that $P$ is strongly irreducible and $M$ is irreducible.  Following \fullref{lemma:boundstab}, we may as well assume that either $M$ contains an essential closed surface of genus less than $\mathrm{genus}(Q)$ or $Q$ is strongly irreducible.  In the former case, the result follows from
Hartshorn's theorem \cite{Ha} (see also \cite{Sc}).  So we are left with the case in which $Q$ is strongly irreducible.

$P$ and $Q$ are not well-separated, for if they were there would be an
incompressible closed surface, not boundary parallel, lying entirely
on one side of $Q$, ie inside a compression body.  As in the proof of \fullref{cor:Heegmain} the result follows
immediately from \fullref{theorem:distance}.
\end{proof}

\section{Appendix:  A quadrilateral Sperner's Lemma}

Here we sketch two different proofs of \fullref{lemma:sperner}.
The first uses algebraic topology and the second is completely
elementary.  

\bigskip

{\bf Algebraic Topology Proof}\qua  Consider $\bdd(I \times I)$ with its standard 
simplicial structure as a $4$--gon.  The dual complex $K$ is also a 
$4$--gon, with a vertex in the interior of each side of $\bdd(I \times 
I)$.  There is a 
natural way to label the vertices of $K$: the vertex in the 
East side of $\bdd(I \times I)$ is labelled $E$ and so on.  

\begin{figure}[ht!]\small
\centering
\labellist
\pinlabel {$\partial(I\times I)$} [b] at 27 169
\pinlabel $K$ [tl] at 45 127
\hair 4pt
\pinlabel $N$ [t] at 80 160
\pinlabel $S$ [b] at 78 19
\pinlabel $E$ [r] at 148 90
\pinlabel $W$ [l] at 10 88
\endlabellist
\includegraphics[width=0.3\textwidth]{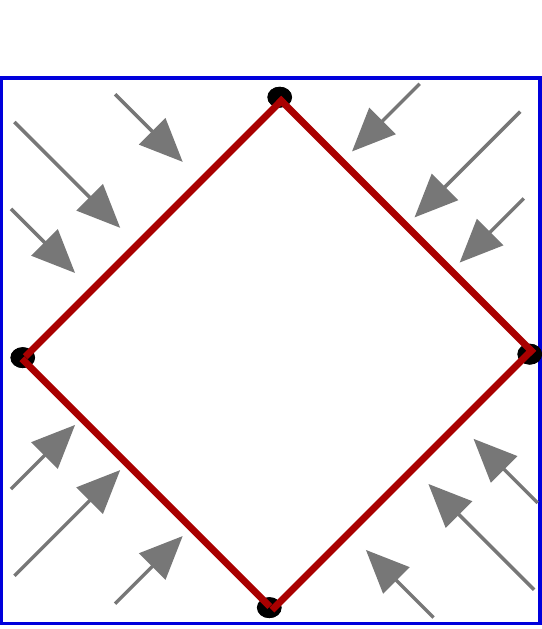}
\caption{} \label{fig:identity}
\end{figure}

{\em Side note:} Of course the identity map $\bdd(I \times I) \map K$ is not a simplicial map; a 
useful way to think of it as a map from one $4$--gon to another is 
shown in \fullref{fig:identity}.

Now consider the subdivision of $I \times I$ given by the labelled
graph $\Lll$.  We will assume that no quadrilateral face of $\Lll$ has
corners containing all $4$ labels and arrive at a contradiction.  This
assumption means that at least two vertices of each quadrilateral have
the same label.  If they are diagonally opposite, add an edge between
them, triangulating the quadrilateral face.  If they are adjacent, add
either diagonal.  The resulting triangulation $\Lll_{+}$ of $I \times
I$ still has the property that no edge has ends labelled $E$ and $W$
or $S$ and $N$.  Map each vertex of $\Lll$ to the vertex in $K$ with 
the corresponding label.  Our labelling
convention guarantees that any edge in $\Lll_{+}$ has ends labelled
the same as an edge in $K$ and, moreover, no $2$--simplex in $\Lll_{+}$
has more than two labels.  It follows that the map of the vertices of
$\Lll$ extends to a simplicial map $\rho: \Lll_{+} \map K$.

Now examine the map $\rho|\bdd( I \times I)$.  Our labelling rule
guarantees that $\rho$ maps the entire East side of $I \times I$ to
the two edges of $K$ that are incident to the vertex $E$ in $K$, ie to
$star(E)$.  Similarly for the other three sides.  In particular, 
$\rho(x)$ is never the antipode of $x$, so $\rho$ is homotopic to the 
identity map.  But it is a classical theorem (equivalent to the 
Brower Fixed Point Theorem) that a map of the circle to itself that is 
homotopic to the identity does not extend over the disk. The 
contradiction completes the proof.~$\sq$\break

\begin{figure}[ht!]
\centering
\labellist
\pinlabel $N$ [b] <2pt,1pt> at 182 337
\pinlabel $S$ [t] at 184 30
\pinlabel $E$ [l] at 339 184
\pinlabel $W$ [r] at 28 186
\endlabellist
\includegraphics[width=0.7\textwidth]{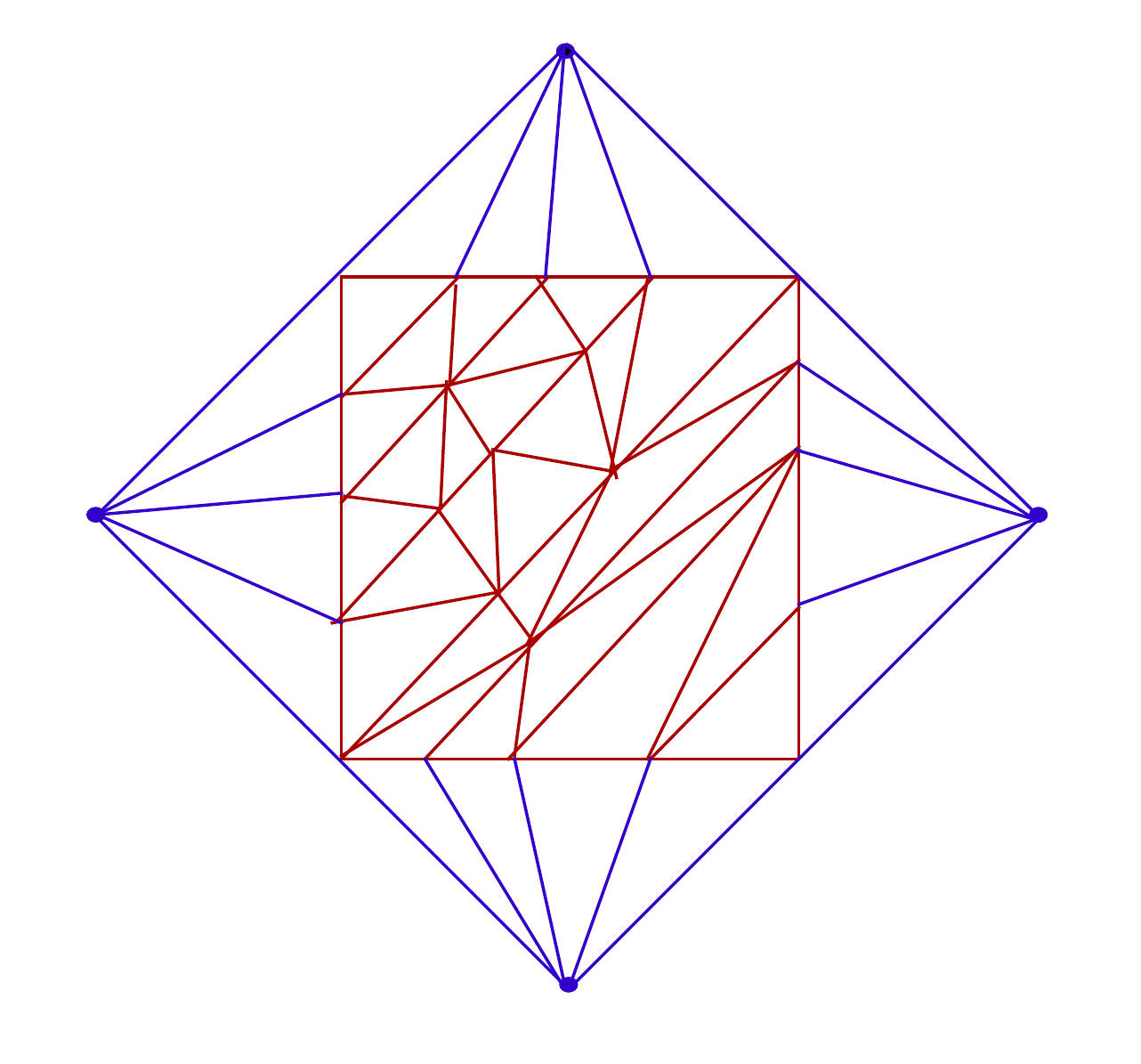}
\caption{} \label{fig:sperner3}
\end{figure}  

\medskip
{\bf Elementary Proof}\qua This proof is modelled on elementary proofs of
Sperner's Lemma (following very helpful comments of Francis Su, cf
also \cite{DPS}).  As above, assume there is no such quadrilateral
face, and again extend $\Lll$ to the triangulation $\Lll_{+}$ as
described in the first proof.  We will give an elementary argument
that there can be no triangulation of the interior of $I \times I$
that satisfies the stated labelling rules on edges and on the sides of
$\bdd (I \times I)$.

\medskip 
Place four extra vertices in the plane outside the square, with one
labelled $N$ directly north of the square, one labelled $E$ directly
east, etc.  Cone each of the new points to the corresponding side of
the square, creating a collection of new triangles on the outside of
the square.  Call the union of these triangles and the original square
a triangulated {\em diamond}.  The diamond has four sides, each with
two edges.  Because no vertex on the North side of the square is
labelled $S$ (and symmetrically for the other four sides) the
labelling on each edge in the diamond still satisfies the given
labelling rule on edges.  Moreover, since the Northwest corner of the
original square must have been labelled either $N$ or $W$, exactly one
of the edges on the Northwest side of the diamond has exactly the
labels $N$ and $W$.  See \fullref{fig:sperner3}.  Symmetrically, of
the eight edges on the boundary of the diamond, exactly one has labels
$S$ and $W$, one has labels $S$ and $E$ and one has labels $N$ and
$E$.  Each of the other four edges has both ends labelled the same (eg
both labelled $N$).

The labelling rule on edges ensures that in each triangle into which
the diamond is divided, at most two labels appear.  Let $t$ be the
number of triangles in the triangulation containing both the labels
$N$ and $W$ and $e$ the number of edges in the triangulation
containing both these labels.  Any edge has two sides in the plane so,
in the plane, there are $2e$ sides of edges labelled exactly $N$ and
$W$.  $2e - 1$ of these sides of edges lie inside the diamond, since,
as we have seen, exactly one edge of the boundary of the diamond
contains the labels $N$ and $W$.  On the other hand, each triangle
that contains labels $N$ and $W$ contains exactly two sides of edges
labelled $N$ and $W$.  Thus $2t = 2e - 1$, a parity contradiction.
\qed

\bibliographystyle{gtart}
\bibliography{link}

\end{document}